\documentclass[11pt,oneside]{amsart}

\usepackage{amssymb,amsmath}
\usepackage{hyperref}
\usepackage{amsthm,amsfonts}
\usepackage{mathrsfs}
\usepackage{hyperref}
\hypersetup{colorlinks,allcolors=blue}
\usepackage{footnote}
\usepackage{textcomp}
\usepackage[english]{babel}
\usepackage{graphicx, color}
\usepackage{tgschola}                     % Dummytext
\usepackage{xargs}                    
\usepackage[pdftex,dvipsnames]{xcolor}  % Coloured text etc.
\usepackage[colorinlistoftodos,prependcaption,textsize=tiny]{todonotes}
\newcommandx{\unsure}[2][1=]{\todo[linecolor=red,backgroundcolor=red!25,bordercolor=red,#1]{#2}}
\newcommandx{\change}[2][1=]{\todo[linecolor=blue,backgroundcolor=blue!25,bordercolor=blue,#1]{#2}}
\newcommandx{\info}[2][1=]{\todo[linecolor=OliveGreen,backgroundcolor=OliveGreen!25,bordercolor=OliveGreen,#1]{#2}}
\newcommandx{\improvement}[2][1=]{\todo[linecolor=Plum,backgroundcolor=Plum!25,bordercolor=Plum,#1]{#2}}
\newcommandx{\thiswillnotshow}[2][1=]{\todo[disable,#1]{#2}}
\numberwithin{equation}{section}
\newtheorem{theorem}{Theorem}[]

\theoremstyle{definition}
\newtheorem{definition}[theorem]{Definition}

\newtheorem{remark}[theorem]{Remark}
\newtheorem{example}[theorem]{Example}

\newtheorem{acknowledgement}{Acknowledgement}

\numberwithin{equation}{section}

\newcommand{\PP} {{\mathbb{P}}}

\usepackage{relsize}
\begin{document}

\title[Hadamard]{Hadamard}
\thanks{}

\author{Iman Bahmani Jafarloo}
%\address{DISMA, Dipartimento di Scienze Matematiche, Politecnico di Torino, Corso Duca degli Abruzzi 24, 10129 Turin, Italy. \newline \hspace*{0.3cm}  Dipartimento di Matematica, Universit\`a degli Studi di Torino, Via Carlo Alberto 10, 10123 Turin, Italy.}
\email{iman.bahmanijafarloo@polito.it}

\keywords{Hadamard Product, Segre Product, Macaulay2}
\subjclass[2020]{}
\thanks{Update: 11 Dec 2020}

\begin{abstract}
This is the \emph{Hadamard} package for Macaulay2 which computes the \emph{Hadamard product} of projective subvarieties.
\end{abstract}

\maketitle

\section*{Introduction}
The \emph{Hadamard product} is a sort of multiplication of matrices and in spite of the usual product it is commutative. The Hadamard product is called as the entry wise or \emph{Schur product} too. The Hadamard product of matrices is a classical topic in linear algebra that is useful in statistics and physics. Also, the product is studied and used in combinatorial and probabilistic problems. Its applications can be found not only in mathematics but also for instance, in cryptography, information theory and data compression such as jpeg format. Around 2010, in \cite{CMS,CTY}, the Hadamard product of matrices was extended to Hadamard product of varieties in the study of the geometry of \emph{binary Boltzmann machine}. The Hadamard product of varieties is also related to \emph{tropical geometry} \cite{BCK, FOW, MS}. Other papers that contributed to the study of Hadamard product of varieties include \cite{BJC,BCFL1,BCFL2,CCFL,CCGV}.

This text is a tutorial paper for the \verb|Hadamard.m2| package of \verb|Macaulay2| \cite{GS} which provides a brief introduction of the Hadamard product of varieties and several examples. The package has been verified for the Macaulay2 development and it can be found at \url{https://github.com/imanbj/Hadamard-m2}.
\section*{Hadamard Product}
\begin{definition}[\cite{CMS}]
Given any two subvarieties $X$ and $Y$ of a projective space $\mathbb{P}^{n}$, their usual Segre product is
$$\mathbb{P}^{N} \times \mathbb{P}^{N} \longrightarrow \mathbb{P}^{N}$$
$$\left(X,Y\right) \mapsto X\times Y,$$
where $[x_0y_0:x_0y_1:\cdots:x_ny_n]\in X \times Y\subset\PP^N $.
We denote by $ z_{ij} $ the coordinates of $ \PP^N $. We define their \textit{Hadamard product} $ X \star Y $ to be the closure of the image of the rational map  $\pi: \mathbb{P}^{N} \dashrightarrow \mathbb{P}^{n}$ 
from the linear space $ \Gamma $ defined by equations $ z_{ii}=0 $, $ i =0,\ldots, n $. Hence, we have that $ X\star Y = \overline{\pi(X\times Y)} $. Note that the closure is the Zariski's closure.
\end{definition}
\begin{definition}
For two points $x=\left[x_{0}: \cdots: x_{n}\right]$ and $y=\left[y_{0}: \cdots: y_{n}\right]$  in $\mathbb{P}^{n}$ with $a_{i} b_{i} \neq 0$ for some $i$, we define the Hadamard product $x \star y$ of $x$ and $y$ as
$$
x \star y=\left[x_{0} y_{0}: x_{1} y_{1}: \cdots: x_{n} y_{n}\right].$$
If $x_{i} y_{i}=0$ for all $i=0, \ldots, n$ then we say $x \star y$ is not defined.
\end{definition}
In this package a \emph{point} is constructed from a list \verb|{}| or an array \verb|[]|, see Example \ref{point}. 
\begin{example}\label{point}
Consider the points $[1 : 2 : 3],\ [-1:2:5],\ [1:1:2],\ [2:2:4] $.
	\begin{verbatim}
	Macaulay2, version 1.14
	i1 : loadPackage"Hadamard"
	i2 : p = point {1,2,3};
	i3 : q = point [-1,2,5];
	i4 : p * q
	o4 = Point{-1, 4, 15}
	i5 : x = point {1,1,2};
	i6 : y = point {2,2,4};
	i7 : x == y
	o7 = true
	i8 : p == q
	o8 = false
	\end{verbatim}
\end{example}
\begin{remark}\label{closure}
	We have that
	$$X \star Y=\overline{\{x \star y: x \in X, y \in Y, x \star y \text { is defined }\}}.$$
Note that the closure of the map $ \pi $ is necessary unless it is defined on all of $ X\times Y $, see \cite[Remark 2.4]{BCK} and Example \ref{exclosure}.
\end{remark}

Let $ X $ and $ Y $ be two varieties in $ \PP^n $. Denote by $ I:=I(X)\subset R=K[x_0,\dots,x_n] $ and $ J:=J(Y)\subset T=K[y_0,\ldots,y_n] $ their corresponding vanishing ideals. We consider 
$$S = \frac{R}{I}\otimes\frac{T}{J}.$$
Let $ W=K[w_0,\ldots,w_n] $. We define the map  $$ \Phi:W \longrightarrow S $$
$$w_0\mapsto x_0y_0$$
$$\vdots$$
$$w_n\mapsto x_ny_n.$$
Hence, we have that $ X\star Y = \ker \Phi$. 

The function \verb|hadamardProduct| is multi functional and works in different ways. It computes the Hadamard products of two ideals, two sets of points, and a list of ideals or  points.  
\begin{example}
Consider the ideals $ I$ and $J $.
\begin{verbatim}
i9  : S = QQ[x,y,z,w];
i10 : I = ideal(x-y+z,z+y+w);
i11 : J = ideal(2*x-y,w+x+y+z);
i12 : hadamardProduct(I,J)           
o12 = ideal(6x^2 - 9x*y + 3y^2 - 2x*z + 2y*z + 2x*w - y*w)
\end{verbatim}
Using \verb|idealOfProjectivePoints|, one can compute the vanishing ideal of a list of points.
\begin{verbatim}
i13 : p = point {1,-1,4,2};
i14 : q = point {1,2,0,-1};
i15 : idealOfProjectivePoints({p*q},S)
o15 = ideal (z, y - w, 2x + w)
\end{verbatim}
This can be computed also from their defining ideals
\begin{verbatim}
i23 : Ip = idealOfProjectivePoints({p},S)
o23 = ideal (z - 2w, 2y + w, 2x - w)
i24 : Iq = idealOfProjectivePoints({q},S)
o24 = ideal (z, y + 2w, x + w)
i25 : hadamardProduct(Ip,Iq)
o25 = ideal (z, y - w, 2x + w)
i26 : o25 == o15
o26 = true
\end{verbatim}
\begin{example} \label{exclosure}
Due to Remark \ref{closure},  consider the curve $X \subset \mathbb{P}^{2}$ defined by $x y-z^{2}=0$ and the point $p=[1: 0: 2]$. The Hadamard product $p \star X$ is the line $y=0 .$ However the point $[0: 0: 1]$ cannot be obtained as $p \star q,$ with $q \in X$.

\begin{verbatim}
i27 : T = QQ[x,y,z];
i28 : IX = ideal(x*y-z^2)
i29 : p = point[1, 0, 2]
i30 : Ip = idealOfProjectivePoints({p},T)
o30 = ideal (y, 2x - z)
i31 : hadamardProduct(Ip,IX)
o31 = ideal y
\end{verbatim}

\end{example}

\end{example}
Given two sets of points $ L $ and $ M $ returns the list of (well-defined) entry-wise multiplication of pairs of points in the Cartesian product $ L\times M $.
\begin{example}
Consider the two sets of points $ L$ and $M $.
\begin{verbatim}
i32 : L = point\{{0,1}, {1,2}};
i33 : M = point\{{1,0}, {2,2}};
i34 : hadamardProduct(L,M)
o34 = {Point{1, 0}, Point{0, 2}, Point{2, 4}}
\end{verbatim}
Note that the map $ \pi $ is not defined on all $ L\times M $ since \verb|point{0,1} * point{1,0} = point{0,0}|. 
\end{example}
The Hadamard products of a list of ideals or points are constructed by using iteratively the binary function \verb|hadamardProduct(Ideal,Ideal)|, or \verb|Point * Point|.

\begin{example}
Let $ K $ be an ideal different from $ I $ and $ J $.
\begin{verbatim}
i35 : K = ideal(2*x-y,w+x+y+z, z-w)
o35 = ideal (2x - y, x + y + z + w, z - w)
i36 : L = {I,J,K};
i37 : hadamardProduct(L)
o37 = ideal(72x^2  - 54x*y + 9y^2  + 16x*z - 8y*z - 16x*w + 4y*w)
i38 : P = point\{{1,2,3},{-1,1,1},{1,1/2,-1/3}};
i39 : hadamardProduct(P)
o39 = Point{-1, 1, -1}	
\end{verbatim}
\end{example}
\begin{definition}
For any projective variety $X$, we may consider its Hadamard square $X^{[2]}=X \star X$ and its higher Hadamard powers $X^{[r]}=X \star X^{[r-1]}$.
\end{definition}
Give a homogeneous ideal $ I $ or a set of points $ L $, the $ r $-th Hadamard power of $ I $ or $ L $ is computed by the function \verb|hadamardPower|. The function computes the $ r $-times Hadamard products of $ I $ or $ L $ to itself.
\begin{example}
Let us recall the ideals $ I$ and $K$.
\begin{verbatim}
i40 : hadamardPower(J,3)
o40 = ideal(8x - y)
i41 : hadamardPower(K,3)
o41 = ideal (z - w, 27y + 64w, 27x + 8w)
i42 : L={point{1,1,1/2},point{1,0,1},point{1,2,4}};
i43 : hadamardPower(L,2)
o43 = {Point{1, 0, 1}, Point{1, 0, 1/2}, Point{1, 0, 4},
       Point{1, 4, 16}, Point{1, 2, 2}, Point{1, 1, 1/4}}
\end{verbatim}
\end{example}
\begin{example}
Here is an another example for the Hadamard power.
	
\begin{verbatim}
i44 : S = QQ[x,y,z];
i45 : X = point\{{1,1,0},{0,1,1},{1,2,-1}};
i46 : A = idealOfProjectivePoints(X,S)
o46 = ideal(3x*z - y*z + z^2, 3x*y - 3y^2 - y*z + 4z^2,
            3x^2 - 3y^2 - 2y*z + 5z^2, y^2z + y*z^2 - 2z^3)
i47 : A2 = hadamardPower(A,2)
o47 : ideal(y^2*z - 18x*z^2+ y*z^2 - 2z, x*y*z - 4x*z^2,
            x^2*z - x*z^2, 2x^3 - 3x^2*y + x*y^2 - 6x*z^2) 
i48 : X2 = hadamardPower(X,2);
i49 : A2 == idealOfProjectivePoints(X2,S)
o49 = true
\end{verbatim}
\end{example}

\begin{theorem}[\cite{MS}]
Let $ L $ and $ M $ be two lines in $ \PP^3 $. Then $ L\star M $ is a quadratic form.
\end{theorem}

\begin{example}[\cite{BCK}[Example 6.1]]
Let $L$ be the line in $\mathbb{P}^{3}$ through points $[2: 3: 5: 7]$ and $[11: 13: 17: 19]$ and $M$ the line through $[23: 29: 31: 37]$ and $[41: 43: 47: 53] .$ Let $x, y, z, w$ be the homogeneous coordinates of $\mathbb{P}^{3}$. We compute $L \star M$ as follows:
\begin{verbatim}
i50 : S=QQ[x,y,z,w];
i51 : L = matrix{{2, 3, 5, 7},{11, 13, 17, 19}};
i52 : IL = ideal flatten entries(matrix{gens S} * gens ker L);
i53 : M = matrix {{23, 29, 31, 37},{41, 43, 47, 53}};
i54 : IM = ideal flatten entries(matrix{ gens S} * gens ker M);
i55 : hadamardProduct(IL,IM)
o55 = ideal(88128x^2 - 89280x*y - 5299632y^2 - 817938x*z + 8896641y*z -
            1481805z^2 - 321510x*w - 1777545y*w - 54250z*w + 116375w^2)
\end{verbatim}

\end{example}
 
\begin{acknowledgement}
The author wishes to thank the organizer of  Macaulay2 workshop, University of Saarlandes, Saarbrücken, Germany, 2019. The author acknowledges that the his traveling and accommodation were fully supported by the workshop. He also thanks Alessandro Oneto for his help in improving the \textit{Hadamard} package especially its documentation.
\end{acknowledgement}

\end{document}